\documentclass[11pt]{amsart}
\usepackage{amssymb}
\usepackage{graphicx}
\usepackage{amsmath, color}
\usepackage{algpseudocode}
\usepackage{psfrag}
\begin{document}

\title[Moment Generating Function for the rectangular half--Gilbert model]{The Moment Generating function for ray lengths in the Half Gilbert Model with Rectangular Cells }

\author{James Burridge*, Richard Cowan**}

\date{\today}

\maketitle

\noindent \begin{center}  *{\scriptsize Department of Mathematics,
University of Portsmouth,
Portsmouth, UK. \textsf{james.burridge@port.ac.uk} }\\
**{\scriptsize School of Mathematics and Statistics, University of
Sydney, NSW, 2006, Australia. \textsf{rcowan@usyd.edu.au} }
\end{center}

\vspace{.6cm}
\begin{abstract}
In the \emph{full} rectangular version of Gilbert's tessellation \cite{eng,mm,bcm}  lines extend either horizontally (with east- and west--growing rays) or vertically (north- and
south--growing rays) from seed points which form a Poisson point
process, each ray stopping when another ray is met. In the
\emph{half} rectangular version \cite{bcm}, east and south growing rays do not
interact with west and north rays. Using techniques developed in our previous paper \cite{bcm}, we derive an exact expression for the moment generating function for the ray length distribution in the half rectangular model.

\end{abstract}

\section{A brief review of stopping sets in the Half Gilbert model}\label{Rev}

Suppose that a stationary Poisson process of intensity $\lambda$
exists in the plane, with seeds marked either $H$ (east growing) with probability $q$ or
$V$ (south growing) with probability $1-q$. These seeds produce rays growing at a constant rate in directions matching their label ($H$ or $V$). Rays end their growth at the instant that the growing tip meets another ray. Our aim is to investigate the terminal length distribution of an east growing test ray, whose seed location we take to be the origin $O$. The seeds that are relevant for the test ray lie in the \emph{unbounded} octant between the lines $y=x$ and $y=0$, with $x\geq 0$: we call this region, the initial \emph{live zone}.

\begin{figure}[h]
    \psfrag{r}{$r$}\psfrag{A}{$S_1$}\psfrag{B}{$D_1$}\psfrag{C}{$S_2$}
        \psfrag{y}{$y$}\psfrag{P}{$O\ $}
\begin{center}
\includegraphics[width=8 cm]{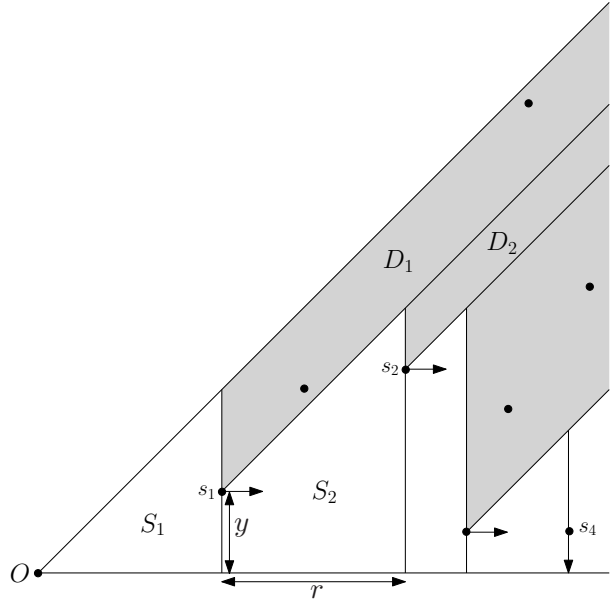}
\caption{\scriptsize \label{trapSS} Trapezoidal stopping sets and
dead zones in the half model. }
\end{center}
\end{figure}

We construct a stopping set process \cite{zuy99,cowan03} by expanding a domain --- an isosceles right angle triangle (see Figure \ref{trapSS}) --- into the
live zone, stopping when it hits the first seed $s_1$ whose
coordinates relative to $O$ are $(x_1, y_1)$. This creates a domain
$S_1$ with area $E_1$ that is exponentially distributed. If $s_1$ is
$V$--type, then it will provide the ray that blocks the test seed;
thus $L=x_1$ and no other seeds need be considered.

Alternatively if $s_1$ is $H$--type, then, instead of growing $S_1$
(retaining its shape as an isosceles right--angle triangle), we remove a part of the live zone: a `dead zone' labelled $D_1$ (see Figure \ref{trapSS}) which
has now become irrelevant. As $S_1 \cup D_1$ has been constructed without drawing upon any information taken from outside $S_1 \cup D_1$, the point process in
the remaining region (the new \emph{live zone}) is still a Poisson
process with unchanged intensity given the information within $S_1
\cup D_1$.

We now grow a trapezium whose left--hand side located at $x=x_1$ has
length $y=y_1$. The trapezium expands until its right--hand side
first hits a seed $s_2$ (in the new live zone). The stopping set
formed is called $S_2$. It has an exponentially distributed area
$E_2$.

We proceed in this way, forming a sequence of stopping sets
(illustrated in Figure \ref{trapSS}) with independent
exponentially--distributed areas. The sequence possesses a recursive
structure, which we can exploit. It is also important to note that
the first $V$--type seed will provide the ray which blocks the test
ray. Without our introduction of dead zones, a complicated algorithm would be required to check if a $V$--type ray actually reaches the path of the test ray.

\vspace{.4cm}\textbf{The recursive structure commencing with a
generic live zone:} Suppose that we begin observing the process when
the live zone has left boundary of height $y$ and when we are about
to construct $S_n$. In Figure \ref{trapSS}, we draw the case $n=2$.
The probability density function for the length, $r$, of $S_n$'s
base, conditional on the height $y$ of its left boundary, follows
from the exponential distribution of $S_n$'s area $E_n$. It is
therefore:
\begin{equation*}
f(r\mid y) = \lambda (r+y)e^{-\frac{\lambda}{2}(r^2+2ry)}.
\end{equation*}
If the stopping seed $s_n$ for set $S_n$ is $V$--type, then its
south ray will be the first to intersect the test ray and the
process ends. Otherwise, another dead zone is created and further
trapezoidal stopping sets are formed until a $V$--type is met.

Let $X$ be the random variable equal to the horizontal distance
covered by stopping sets until the process comes to an end. The
density function of $X$, conditional on $y$ will be:
\begin{equation*}
\label{conDensity}
g(x|y) = (1-q) \lambda (x+y)e^{-\frac{\lambda}{2}(x^2+2xy)} + q\lambda \int_0^\infty e^{-\frac{\lambda}{2}(r^2+2ry)}\left[\int_0^{r+y} g(x-r|u) du \right] dr.
\end{equation*}
where $g(x|y)=0$ if $x<0$. The first term accounts for the case where the first seed is V--type, and the second term for the case where it is H--type and the process is
effectively re--started with a different boundary condition having
already covered some horizontal distance. Note that the ray length
probability density function is $g(x\mid 0)$.

\section{The Moment Generating Function}

We define the moment generating function of the conditional density $g(x \mid y)$:
\begin{equation*}
M_t(y)=\int_0^\infty e^{tx} g(x|y) dx,
\end{equation*}
which, from equation (\ref{conDensity}), satisfies the integral equation:
\begin{multline}
\label{M1}
M_t(y) = (1-q) \int_0^\infty (x+y) e^{-\frac{1}{2}(x^2+2xy)} dx \\ + q \int_0^\infty e^{-\frac{1}{2}(r^2+2ry)} \left[ \int_0^{r+y} \left( \int_0^\infty e^{xt} g(x-r \mid u) dx \right) du \right] dr.
\end{multline}
To solve this equation we use the following method. We begin by converting equation (\ref{M1}) into a differential equation whose general solution is a combination of Kummer confluent hypergeometric functions. Retaining only those parts of the solution which have the correct asymptotic behaviour, it remains to deduce the form of a single arbitrary function of $t$. This function is found by substituting the solution back into the integral equation evaluated at $y=0$.

Commencing with the first part of our method, we note that the second term of (\ref{M1}) contains the integral:
\begin{align*}
\int_0^\infty e^{xt} g(x-r \mid u) dx &= \int_0^\infty e^{(z+r)t} g(z \mid u) dz \\
&= e^{rt} M_t(u).
\end{align*}
We may now simplify the second term using integration by parts:
\begin{multline*}
\int_0^\infty e^{-\frac{1}{2}(r^2+2r(y-t)} \left[ \int_0^{r+y} M_t(u) du \right] dr = \\ e^{\frac{1}{2}(y-t)^2} \sqrt{\frac{\pi}{2}} \left[ \textrm{erfc} \left(\frac{y-t}{\sqrt{2}} \right) \int_0^y M_t(u) du + \int_y^\infty \textrm{erfc} \left(\frac{u-t}{\sqrt{2}} M_t(u) du \right) \right].
\end{multline*}
Evaluating the integral in the first term of (\ref{M1}), and then combining both terms in together we obtain:
\begin{multline*}
\label{M2}
M_t(y) = (1-q)\left[1 + \sqrt{\frac{\pi }{2}} t e^{\frac{1}{2} (t-y)^2}  \text{erfc} \left(\frac{y-t}{\sqrt{2}}\right)\right] \\ + qe^{\frac{1}{2}(y-t)^2} \sqrt{\frac{\pi}{2}} \left[ \textrm{erfc} \left(\frac{y-t}{\sqrt{2}} \right) \int_0^y M_t(u) du + \int_y^\infty \textrm{erfc} \left(\frac{u-t}{\sqrt{2}}\right) M_t(u) du  \right].
\end{multline*}
This may be reduced to the differential equation:
\begin{equation*}
\frac{d^2M_t(y)}{dy^2} - (y-t) \frac{dM_t(y)}{dy} - (1-q) M_t(y) = -(1-q),
\end{equation*}
which has the general solution:
\begin{equation}
\label{genSol}
M_t(y) = 1 + c(t) H_{q-1}\left(\frac{y-t}{\sqrt{2}}\right) + d(t) \,
   _1F_1\left(\frac{1-q}{2};\frac{1}{2};\left(\frac{y}{\sqrt{2}}-\frac{t
   }{\sqrt{2}}\right)^2\right),
\end{equation}
where $c(\cdot)$ and $d(\cdot)$ are arbitrary functions of $t$. The last term of expression (\ref{genSol}) is an example of a Kummer confluent hypergeometric function. Because all the $t$ derivatives of this function, evaluated at $t=0$, diverge as $y \rightarrow \infty $, we may discard it from the solution on the grounds that for all $x>0$, $\lim_{y \rightarrow \infty} g(x \mid y) =0 $. The second term in the solution, a generalisation of the Hermite polynomials, $H_n(x)$, to non-integer $n$, is in fact a combination of Kummer functions with the correct asymptotic behaviour \cite{bcm}. For our purposes, we will need the following integral representation of the Hermite function:
\begin{equation}
\label{HIntRep}
H_v(z) = \frac{2^{v+1}}{\sqrt{\pi}} e^{z^2} \int_0^\infty e^{-u^2}u^v \cos \left( 2 z u - \frac{\pi v}{2} \right) du.
\end{equation}

Now that we have have the correct general form of the moment generating function, it remains to evaluate the constant $c(t)$. To do so, we substitute the solution (\ref{genSol}) with $d(t)=0$, into the original integral equation (\ref{M1}) evaluated at $y=0$. We find that:
\begin{equation*}
\label{M3}
 c(t)=
 \frac{ t\  \text{erfc} \left(\frac{-t}{\sqrt{2}}\right)}{\left[\sqrt{\frac{2}{\pi}}e^{-\frac{1}{2}t^2} H_{q-1}\left(\frac{-t}{\sqrt{2}}\right)-q  \int_0^\infty \textrm{erfc} \left(\frac{u-t}{\sqrt{2}}\right) H_{q-1}\left( \frac{u-t}{\sqrt{2}}\right) du \right]}.
\end{equation*}
In order to find an explicit expression for $c(t)$ we need to evaluate the integral:
\begin{equation*}
J(t,q):=\int_0^\infty \textrm{erfc} \left(\frac{u-t}{\sqrt{2}}\right) H_{q-1}\left( \frac{u-t}{\sqrt{2}}\right) du.
\end{equation*}
To do this, we first make the change of variable $z=(u-t)/\sqrt{2}$, and then split the integration range into two, giving:
\begin{align*}
J(t,q) &= \sqrt{2} \int_{-\frac{t}{\sqrt{2}}}^\infty   \textrm{erfc} (z) H_{q-1}(z) dz \\
 &= \sqrt{2} \left[ \int_{-\frac{t}{\sqrt{2}}}^0  \textrm{erfc} (z) H_{q-1}(z) dz + \frac{\Gamma(-q/2)}{4\Gamma(1-q)}-\frac{2^q}{q^2 \Gamma(-q/2)} \right].
\end{align*}
Now consider the integral:
\begin{equation*}
K(t,q) := \int_{-\frac{t}{\sqrt{2}}}^0  \textrm{erfc} (z) H_{q-1}(z) dz.
\end{equation*}
Integration by parts gives:
\begin{equation*}
K(t,q) = \frac{1}{2q} \left[\textrm{erfc} (z) H_q(z)\right]_{-\frac{t}{\sqrt{2}}}^0 + \frac{1}{q\sqrt{\pi}} \int_{-\frac{t}{\sqrt{2}}}^0  e^{-z^2} H_q(z) dz.
\end{equation*}
In order to evaluate the second term, we make use of the integral representation of the Hermite function (\ref{HIntRep}), finding that:
\begin{align*}
\int_{-\frac{t}{\sqrt{2}}}^0  e^{-z^2} H_q(z) dz &=  \frac{2^{q+1}}{\sqrt{\pi}} \int_0^\infty e^{-u^2}u^q\left[ \int_{-\frac{t}{\sqrt{2}}}^0  \cos \left( 2 z u - \frac{\pi v}{2} \right)dz \right] du \\
&= \frac{2^{q+1}}{\sqrt{\pi}} \int_0^\infty e^{-u^2}u^{q-1} \sin \left(\frac{t u}{\sqrt{2}}\right) \cos \left(\frac{1}{2}
   \left(\pi  q+\sqrt{2} t u\right)\right) du \\
&= \frac{2^{q-\frac{1}{2}}}{\sqrt{\pi }} t \cos \left(\frac{\pi  q}{2}\right) \Gamma
   \left(\frac{q+1}{2}\right) \,_1F_1\left(\frac{q+1}{2};\frac{3}{2};-\frac{t^2}{2}\right)\\
  & + \frac{2^{q-1}}{\sqrt{\pi }} \sin
   \left(\frac{\pi  q}{2}\right) \Gamma \left(\frac{q}{2}\right)
   \left(\,
   _1F_1\left(\frac{q}{2};\frac{1}{2};-\frac{t^2}{2}\right)-1\right).
\end{align*}
The function $K(t,q)$, in full, is therefore:
\begin{multline*}
K(t,q) = \frac{1}{2q} \left[\frac{2^q \sqrt{\pi}}{\Gamma\left(\frac{1-q}{2}\right)}-\textrm{erfc}\left(-\frac{t}{\sqrt{2}}\right)
H_q\left(-\frac{t}{\sqrt{2}}\right)\right]\\
+\frac{2^q}{2q\pi} \left[ \sqrt{2} t \cos \left(\frac{\pi q}{2}\right) \Gamma\left(\frac{q+1}{2}\right)\,_1F_1\left(\frac{q+1}{2};\frac{3}{2};-\frac{t^2}{2}\right)
 \right.\\
 \left. + \sin \left(\frac{\pi  q}{2}\right) \Gamma \left(\frac{q}{2}\right) \left(\,
   _1F_1\left(\frac{q}{2};\frac{1}{2};-\frac{t^2}{2}\right)-1\right)\right],
\end{multline*}
and $J(t,q)$ is:
\begin{equation*}
J(t,q) = \sqrt{2} \left[ K(t,q) + \frac{\Gamma(-q/2)}{4\Gamma(1-q)}-\frac{2^q}{q^2 \Gamma(-q/2)} \right].
\end{equation*}
Putting this all together we have:
\begin{equation*}
M_t(y) = 1 +\frac{ t\  \textrm{erfc} \left(\frac{-t}{\sqrt{2}}\right)H_{q-1}\left(\frac{y-t}{\sqrt{2}}\right)}{\left[\sqrt{\frac{2}{\pi}}e^{-\frac{1}{2}t^2} H_{q-1}\left(\frac{-t}{\sqrt{2}}\right)-q J(t,q) \right]}.
\end{equation*}
Setting $y=0$ we have an exact expression for the moment generating function for ray lengths in the half--rectangular Gilbert model. It is worth noting that our analysis has not lost any generality by considering an east moving ray. To obtain the moments of a south moving ray we make the switch $q \leftrightarrow 1-q$. Using this expression we find that the first four moments of terminal ray length are:
\begin{align*}
\mu_1 &= \frac{\Gamma \left(\frac{1}{2}-\frac{q}{2}\right)}{\sqrt{2} \Gamma
   \left(1-\frac{q}{2}\right)} \\
\mu_2 &= \frac{q \Gamma \left(\frac{1}{2}-\frac{q}{2}\right)^2}{\Gamma
   \left(1-\frac{q}{2}\right)^2}+2 \\
\mu_3 &= \frac{3 \Gamma \left(\frac{1}{2}-\frac{q}{2}\right)}{\sqrt{2} \Gamma
   \left(1-\frac{q}{2}\right)}  \left[ 1 + 2q + \left(\frac{q \Gamma\left(\frac{1}{2}-\frac{q}{2}\right)}{\Gamma
   \left(1-\frac{q}{2}\right)}\right)^2\right] \\
\mu_4 &= 8 \left[1+ q+q (1+  2q)\left(\frac{ \Gamma\left(\frac{1}{2}-\frac{q}{2}\right)}{\Gamma
   \left(1-\frac{q}{2}\right)}\right)^2+\frac{12}{q} \left(\frac{ \Gamma\left(\frac{1}{2}-\frac{q}{2}\right)}{\Gamma
   \left(1-\frac{q}{2}\right)}\right)^4 \right]
\end{align*}
These have a pleasing form, all being polynomial is the same ratio of Gamma functions. We note that the expression $\mu_1$ represents a considerable simplification of our previous formula for this moment \cite{bcm}. To verify that our analysis is correct, let us compare these exact expressions to the approximate moments calculated using the first 200 coefficients from the Cowan--Ma recurrence when $q=\tfrac{2}{5}$. Table \ref{test} shows the results of these calculations, which indicate that our exact expressions are correct.

\begin{table}

\begin{tabular}{|c|c|c|}
  \hline
  Moment & Recurrence & Exact \\
    \hline
 1 & 1.81696 & 1.81696 \\
 2 & 4.64107 & 4.64107 \\
 3 & 15.57 & 15.5701 \\
 4 & 65.9719 & 65.9721 \\
 5 & 342.236 & 342.243 \\
  \hline
\end{tabular}
\caption{\label{test}The first five moments of terminal ray length when $q=\tfrac{2}{5}$ computed using the first 200 coefficients from the Cowan--Maa recurrence relation, compared to the exact results from the moment generating function.}
\end{table}

We conclude by noting that when $q=\tfrac{1}{2}$ the moment generating function has a particularly simple form:
\begin{equation*}
M(t) = 1-\frac{4 \sqrt{2 \pi } t H_{-\frac{1}{2}}\left(-\frac{t}{\sqrt{2}}\right)}{\Gamma
   \left(-\frac{1}{4}\right) \,
   _1F_1\left(-\frac{1}{4};\frac{1}{2};\frac{t^2}{2}\right)+\sqrt{2} t \Gamma
   \left(\frac{1}{4}\right) \,
   _1F_1\left(\frac{1}{4};\frac{3}{2};\frac{t^2}{2}\right)}.
\end{equation*}

\section{Concluding comment}

We have determined the exact form of the moment generating function for terminal ray length in the rectangular half--Gilbert model. To the authors' knowledge, this model remains the only Gilbert--style model for which analytical results exist.

\end{document}